\documentclass[12pt,a4paper]{article}
\usepackage{amsmath,amsthm,amsfonts,amssymb}
\def\ni{\noindent}

\theoremstyle{definition}

\newcommand{\eps}{\epsilon}
 
\newcommand{\Dt}{\Delta} 
\newcommand{\Z}{\mathbb{Z}}
\newcommand{\N}{\mathbb{N}}

\def\ot{\otimes}
\def\wt{\widetilde}
\def\lr{\longrightarrow}
\def\b{\bar}
\def\Ker{{\rm Ker}}
\def\l{\lambda}
\def\eps{\epsilon}
\def\al{\alpha}
\begin{document}
	\begin{center}
		\large{\bf The ring of twisted tilting modules of quantum $GL_n$}\\
		\vspace{0.25cm}
		\small{ M. S. Datt and Rubayya}\\
			{\it  School of Mathematics and Statistics\\
				University of Hyderabad,   Hyderabad, India-500046.\\
				e-mail: msdatt@uohyd.ac.in}
		
	\end{center}

	Abstract: {  Let $G$ be the quantum $GL_n$ over a field of characteristic $p\neq 0$. In this paper we define the ring of twisted tilting modules of $G$.  We give generators and relations for  the ring of twisted tilting modules of quantum $GL_2(k)$. We also show that this is a reduced ring.}\\
	
	\ni {\it Key words: Quantum groups, Tilting modules, Twisted tilting modules, Representation ring.}
	
	\ni \ni {\it 2010 Mathematics Subject Classification 20G05.}

	\section{Introduction}
	Let $k$ be a field of characteristic $p \neq 0$ and $q$ be a primitive $\ell^{th}$ root of unity. Let $G$ be the quantum $GL_n$ over $k$ at $q$ as described in \cite{Dip}. We know that the coordinate ring $k[G]$ is a Hopf algebra. Let $c_{ij}, 1 \leq i,j \leq n$, be the coordinate functions of $G$. Let $T$ be a maximal torus and $B$ be a Borel subgroup of $G$ containing $T$. Therefore the Hopf ideal corresponding to $T$ is given by $\{c_{ij}~|~i\neq j\}$ and the Hopf ideal corresponding to $B$ is given by $\{c_{ij}~ | ~1\leq i<j\leq n\}$.\\
	 	 
	 Let $X(T)={\mathbb{Z}}^n$. The group $X(T)$ is called the character group of $T$. The integral group ring ${\mathbb{Z}}X(T)$ with basis $\{e(\lambda)~|~\lambda\in X(T)\}$ is a group algebra with multiplication given by the rule $e(\l)e(\mu)=e(\l+\mu)$ for $\l, \mu \in X(T)$. Given $\lambda=(\lambda_1,\lambda _2\dots,\lambda_n)\in X(T)$, we have a 1-dimensional $T$-module $k_{\lambda}$ with the comodule structure map sending $a$ to $a\otimes c_{11}^{\lambda_1}c_{22}^{\lambda_2}\dots c_{nn}^{\lambda_n} \in k_\lambda \otimes k[T]$. We know that the set $\{k_{\lambda}~|~\lambda\in X(T)\}$ is a complete set of mutually non-isomorphic irreducible $T$-modules. For any quantum group $H$, we denote the set of all finite dimensional $H$-modules by ${\rm mod}(H)$. Given $V \in {\rm mod}(T)$ and $\lambda \in X(T)$, let $V_\lambda$ denote the sum of all submodules of $V$ isomorphic to $k_\lambda$. We call the submodule $V_\lambda$ the $\lambda$-weight subspace of $V$. Any module $V \in {\rm mod}(G)$ is completely reducible as a $T$-module i.e., $V=\oplus _{\l\in X(T)}V_{\l}$. The formal character of $V$ is given by ${\rm ch}V=\sum _{\l\in X(T)}{\rm dim}(V_{\l})e(\l)$. \\
	
	 Let $\eps_i=(0,0,\dots,1,0,\dots,0)\in {\mathbb{Z}}^n$, for $1\leq i\leq n$, where 1 is in the $i^{th}$ position. The set $\Phi=\{\eps_i-\eps _j~|~1\leq i,j\leq n, i\neq j\}$ is called the set of all roots of $G$. We say that a root $\eps _i-\eps _j$ is a positive root if $i<j$ and is a negative root if $i>j$. Let $\Phi^+$ and $\Phi ^-$ be the set of all positive and negative roots respectively. The subset $S=\{\eps_1-\eps_2,\eps _2-\eps_3,\dots ,\eps _{n-1}-\eps_n\}$ of $\Phi$ is called the set of all simple roots. We have a ${\mathbb{Z}}$-bilinear form $(~~,~~)$ on $X(T)$ such that $(\eps _i,\eps _j)=\delta _{ij}$, for $1\leq i,j\leq n$. We define a partial order on $X(T)$ as follows: for $\l ,\mu \in X(T)$, let $\l\leq \mu$ if $\mu-\l$ is a sum of positive roots. Now for each $\l\in X(T)$, we have induced $G$-modules $\nabla (\l)={\rm Ind}_B^Gk_{\l}$. Let $X^+(T)=\{\l\in X(T)~|~\nabla(\l)\neq 0\}$ and it is called the set of all dominant weights. For $\l\in X^+(T)$, the module $\nabla(\l)$ is called a costandard module. For $\l\in X^+(T)$, we have a Weyl or standard module $\Dt (\l)=(\nabla(-w_0\l))^*$ where $w_0$ is the element of $S_n$ with maximum length. For $V\in {\rm mod}(G)$, a filtration $(0)=V_0\subset V_1\subset V_2\subset \cdots \subset V_n=V$ of $G$-submodules of $V$  is called a costandard filtration (or a good filtration) if $V_i/V_{i-1}\cong \nabla (\l_i)$ for some $\l_i\in X^+(T)$ and for all $1\leq i\leq n$. For a module $V\in {\rm mod}(G)$, we write $V\in {\cal F}(\nabla)$ if $V$ admits a good filtration. Similarly for $V\in {\rm mod}(G)$, a filtration $(0)=V_0\subset V_1\subset V_2\subset \cdots \subset V_n=V$ of $G$-submodules of $V$  is called a standard filtration (or a Weyl filtration) if $V_i/V_{i+1}\cong \Dt (\l_i)$ for some $\l_i\in X^+(T)$ and for all $1\leq i\leq n$. We write $V\in {\cal F}(\Delta)$ if $V$ admits a standard filtration. A module $V$ is said to be a tilting module if $V$ admits both costandard and standard filtration. By [\cite{Dongln},p.no. 262], for each $\l\in X^+(T)$ there exists a unique indecomposable tilting module $T(\l)$ of maximum weight $\l$. Also the set $\{T(\l)~|~ \l \in X^+(T)\}$ forms a complete set of mutually non-isomorphic indecomposable tilting modules. Any tilting $G$-module is a direct sum of copies of $T(\l)$'s. \\
	
	We denote the general linear group $GL_n(k)$ over the field $k$ by $\bar{G}$. Let $x_{ij}$ be the coordinate functions of $\bar{G}$ and $d$ be the determinant (ordinary). Let $\bar{T}$ be a maximal torus and $\bar{B}$ be a Borel subgroup containing $\bar{T}$. Let $X(\bar{T})$ be the group of all characters and $\Z X(\bar{T})$ be the integral group ring of $X(\bar{T})$. The coordinate ring $k[\bar{G}]$ of $\bar{G}$ is a Hopf algebra. We fix a Hopf $\mathbb{F}_p$-form of $k[\bar{G}]$ and let $\bar{F}$ be the corresponding Frobenius morphism. We write $\bar{F}^i$ for the $i^{th}$ iterated Frobenius morphism. If $V$ is a $\bar{G}$-module, then we denote the induced module via $\bar{F}$ by $V^{\bar{F}}$. Let $X^+(\bar{T})$ be the set of all dominant weights of $\bar{G}$. Then $\{\bar{L}(\l) ~|~ \l \in X^+(\bar{T}) \}$ be the complete set of mutually inequivalent simple $\bar{G}$-modules. By \cite{Donag}, we also have a complete set $\{\bar{T}(\l) ~|~ \l \in X^+(\bar{T}) \}$ of mutually non-isomorphic indecomposable tilting $\bar{G}$-modules. Any tilting $\bar{G}$-module is a direct sum of $\bar{T}(\l), \l \in X^+(\bar{T})$ and we call them as partial tilting modules. \\
	
	Let $F : G \longrightarrow \bar{G}$ be the morphism whose comorphism takes $x_{ij}$ to $c_{ij}^\ell$ where $1 \leq i,j \leq n$. Thus any $\bar{G}$-module $V$ is also a $G$-module via $F$ and we denote this $G$-module by $V^F$. The morphism $F$ is called the quantum Frobenius morphism. Let $X_1=\{(\l_1, \l_2, \dots,\l_n) \in X^+(T) ~|~ 0 \leq \l_1-\l_2, \l_2-\l_3, \dots, \l_{n-1}-\l_n \leq \ell-1\}$. By the Steinberg's tensor product theorem [\cite{Donqsc},p.no.65]: for $\l \in X_1$ and $\mu \in X^+(\bar{T})$, we have $L(\l + \ell \mu) \cong L(\l)\otimes\bar{L}(\mu)^F$. It has been shown in \cite{Dattsm} that the tensor product of two simple quantum $GL_2(k)$-modules can be expressed as a finite direct sum of indecomposable quantum twisted tilting modules of the form $T_{-1} \otimes (\bar{T}_0 \otimes \bar{T}_1^{\bar{F}} \otimes \cdots \bar{T}_r^{\bar{F}^r})^F$, where $T_{-1}$ is an indecomposable tilting $G$-module and $\bar{T}_i$'s are indecomposable tilting $\bar{G}$-modules, for all $i \geq 0$.\\
	
	Let $\rm {rep}(G)$ denote the set of all isomorphism classes of $\rm {mod}(G)$. For $V \in \rm {mod}(G)$, we denote the class in $\rm {rep}(G)$ through $V$ by $[V]$ . The set $\rm {rep}(G)$ is a ring with respect to the binary operations $[V]+[W]=[V \oplus W]$ and $[V][W]=[V \otimes W]$, for $V,W \in \rm {mod}(G)$. We call the ring $\rm {rep}(G)$ the representation ring of $G$. Let $\mathcal{B}=\{V_i ~|~ i \in I\} \subset \rm {mod}(G)$ be a complete set of mutually non-isomorphic indecomposable $G$-modules. Then $\{[V] ~|~ V \in \mathcal{B}\}$ is a basis for $\rm {rep}(G)$. Similarly we also have the Grothendieck ring $\rm {Grot} G$ of $\rm {mod}(G)$. For $V \in \rm {mod}(G)$, we denote the class in $\rm {Grot} G$ through $V$ by $\langle V \rangle$. Let $V \in \rm {mod}(G)$, we have $\langle V \rangle = \Sigma_{i} a_i \langle L_i \rangle$ where $L_i$'s are simple $G$-modules and $a_i$'s are their composition multiplicities in $V$. For $V,W \in \rm {mod}(G)$, define the multiplication $\langle V \rangle \langle W \rangle = \langle V\otimes W \rangle$. Then the ring $\rm {Grot} G$ is a free $\mathbb{Z}$-module with a basis $ \{ \langle L(\l) \rangle ~|~ \l \in X^+(T) \} $. By [\cite{Donqsc},69,14(iii)], we know that, if $V,W$ are tilting modules, then $V \otimes W$ and $V \oplus W$ are also tilting modules. Hence the span of classes of tilting modules forms a subring of $\rm {rep}(G)$ and we denote it by $\rm {rep_{tilt}(G)}$. Let $\rm {ch}$ denote the map from $\rm {rep} (G)$ to $(\mathbb{Z}X(T))^W$, which takes the class $[V]$ to $\rm {ch}V$. This induces an isomorphism from $\rm {Grot} G$ to $(\mathbb{Z}X(T))^W$. By [\cite{Donag},(1.3)] we have an isomorphism between $\rm {rep_{tilt}(G)}$ and $(\mathbb{Z}X(T))^W$. \\
	
	Let $A(G)$ be the subring of $\rm {rep} (G)$ generated by all tilting $G$-modules and quantum Frobenius twist of twisted tilting modules of $\bar{G}$. In \cite{Doty}, Doty and Henke has proved that the tensor product $L \otimes  L'$ of irreducible rational modules $L,L'$ of $SL_2(k)$ can be expressed as a direct sum of $GL_2(k)$-modules of the form $\bar{T}_0 \otimes \bar{T}_1^{\bar{F}} \otimes \dots \otimes \bar{T}_r^{\bar{F}^r}$, where each $\bar{T}_i$ is an indecomposable tilting module of $GL_2(k)$. In a paper by \cite{Dattag}, it has been shown that the ring of twisted tilting modules of $SL_2(k)$ is a reduced ring. In \cite{Dattsm}, it has been proved that the tensor product of any two simple quantum $GL_2(k)$-modules can be expressed as a finite direct sum of twisted tilting modules. \\
	
	In this paper, we give generators and relations for the ring $A(G)$. We show that the ring of twisted tilting modules for quantum $GL_2(k)$ is a reduced ring.
	
	\section{The ring of twisted tilting modules for quantum $GL_n$}
	 Let $R$ be the free $k$-algebra generated by the elements \{$X_{ij} ~|~ 1\leq i,j\leq n$\}. Then $R$ with comultiplication $\delta: R \longrightarrow R \otimes R, \delta (X_{ij})=\Sigma_{r=1}^n X_{ir} \otimes X_{rj}$ and counit $\epsilon: R \longrightarrow k, \epsilon(X_{ij})=\delta_{ij}$ is a bialgebra. Let $I$ be the ideal of $R$ generated by elements of the form
	\begin{center}
		$X_{ir}X_{is} - X_{is}X_{ir}$ \quad for all $i,r,s$, \\
		$X_{jr}X_{is} - qX_{is}X_{jr}$ \quad for all $i<j$ and $r \leq s$,\\
		$X_{jr}X_{is} - X_{is}X_{jr} - (q-1)X_{ir}X_{js}$ \quad for all $i<j$ and $r>s$,
	\end{center} 
	where $i \leq i,j,r,s \leq n$. By assigning degree as $\deg X_{ij}=1 \ (\text{for} \ 1 \leq i,j \leq n)$, the ideal $I$ is a homogeneous ideal. Let $A_q(n)=R/I$ and $c_{ij}=X_{ij}+I$ for $1 \leq i,j \leq n$. We denote the quantum determinant $\sum _{\sigma\in S_n} sgn(\sigma)c_{1,1\sigma}c_{2,2\sigma}\cdots c_{n,n\sigma}$ by $d_q$. The localisation of $A_q(n)$ with respect to $d_q$ is a Hopf algebra $k[G]$. We denote the corresponding quantum general linear group by $G$. The Hopf algebra $k[G]$ is also called the coordinate ring of $G$. Let $E$ be an $n$-dimensional vector space over $k$ with standard basis $\{e_1,e_2, \dots e_n\}$. The map $e_i \mapsto \sum_{j=1}^{n}e_j \otimes c_{ji}$ from $E$ to $E \otimes k[G]$ defines a left $G$-module structure on $E$. We call the left $G$-module $E$, the natural left $G$-module.  \\ 
	
	Let $\N_0 = \N \cup \{0\}$. Let $\psi^+=\{(i,j) ~|~ 1 \leq j<i \leq n \}$ and $\psi^-=\{(i,j) ~|~ 1 \leq i<j \leq n \}$. Then the $k$-span $I_{\psi^+}$ of all $c_{11}^{a_{11}}c_{12}^{a_{12}}\cdots c_{nn}^{a_{nn}}$, where $(a_{11}, a_{12}, \dots, a_{nn}) \in \N_0^{n^2}$ and $a_{ij} \neq 0$ for some $(i,j) \in \psi^+$, is a Hopf ideal of $k[G]$. The quantum group $B^+$ corresponding to the Hopf algebra $H/I_{\psi^+}$ is called the positive Borel subgroup of $G$. Similarly the quantum subgroup corresponding to the Hopf algebra $k[G]/I_{\psi^-}$ is called the negative Borel subgroup of $G$. Let $\psi_0=\{(i,j) ~|~  i\neq j \}$. Let $I_{\psi_0}$ be the $k$-span of $c_{11}^{a_{11}}c_{12}^{a_{12}}\cdots c_{nn}^{a_{nn}}$, where  $(a_{11}, a_{12}, \dots, a_{nn}) \in \N_0^{n^2}$ and $a_{ij} \neq 0$ for some $(i,j) \in \psi_0$. Then the quantum subgroup $T$ of $G$ corresponding to the Hopf algebra $k[G]/I_{\psi_0}$ is called a torus of the group $G$. Let $G_1$ be the subgroup of $G$ whose defining ideal is the ideal of $k[G]$ generated by $c_{ij}^\ell - \delta_{ij}, 1 \leq i,j \leq n$. Let $B_1=G_1 \cap B, B_1^+=G_1 \cap B^+$ and $T_1=G_1 \cap T$. We let $X_1= \{ \l = (\l_1, \l_2, \dots \l_n) \in X^+(T) ~|~ 0 \leq \l_1-\l_2, \l_2-\l_3,\dots,\l_{n-1}-\l_n, \l_n \leq \ell-1\}$. Then $\{L(\l)_{|G_1} ~|~ \l \in X_1 \}$ is a complete set of mutually non-isomorphic simple $G_1$-modules. \\
	
	Let $\rho=(n,n-1,\dots,1) \in X^+(T)$. We have $L((\ell-1)\rho )= \nabla((\ell-1)\rho) = T((\ell-1)\rho)$ and $\bar{L}((\rho-1)\rho)=\bar{\nabla}((\rho-1)\rho)$. We denote $\nabla((\ell-1)\rho)$ by $\rm {St}$ and $\bar{\nabla}((\rho-1)\rho)$ by $  {\rm \bar{St}}$ and we call them quantum Steinberg and Steinberg (ordinary) modules respectively. More generally, we have the ordinary $r^{th}$ Steinberg module $ {\rm \bar{St}} _r$, for $r \geq 1$. \\
	
	Let $X^+(\bar{T})_\infty$ be the set of all sequences $(\l)= (\l(0),\l(1), \dots)$ of dominant weights of $\bar{G}$ with $\l(i)=0$ for sufficiently large $i$. We let $X^+(T)_\infty =X^+(T) \times X^+(\bar{T})_\infty$. Let $\N_{-1} =\N \cup \{0,-1\}$. For each $n \in \N_{-1}$, we let $X^+(T)_n= \{(\l)=(\l(-1),\l(0),\dots)\in X^+(T)_\infty ~|~ \l(i)=0 \ {\rm for\ all} \ i > n \}$. \\
	
	Let $\phi: {\rm rep} (\bar{G})\longrightarrow {\rm rep}(G)$ be the endomorphism which sends the class $[V]$ to $[V^F]$. We denote $\phi(a)$ by $a^\phi$, for $a \in {\rm rep}(\bar{G})$. Similarly let $\bar{\phi}:{\rm rep}(\bar{G})\longrightarrow {\rm rep}(\bar{G})$ be the endomorphism which sends $[V]$ to $[V^{\bar{F}}]$. As earlier for $\phi$, we denote $\bar{\phi}(a)$ by $a^{\bar{\phi}}$, for $a \in {\rm rep}(\bar{G})$. We write $\bar{a}^{\psi^{r+1}}$ for $(\bar{a}^{\bar{F}^r})^F$, where $r \geq 0$. The ring ${\rm rep}(G)$ has a $\Z$-basis $\{[T(\l)] ~|~ \l \in X^+(T) \}$, by [\cite{Donqsc}, 3.3]. Let us denote the class $[T(\l)]$ by $t(\l)$. Similarly for $\l \in X^+(\bar{T})$, we denote the class $[\bar{T}(\l)]$ by $\bar{t}(\l)$. For $(\l) \in X^+(T)_\infty$, we let $T((\l))= T(\l(-1)) \otimes (\bar{T}(\l(0)) \otimes (\bar{T}(\l(1))^{\bar{F}} \otimes \bar{T}(\l(2))^{\bar{F}^2} \otimes)^F$. Therefore, we have $t((\l))=t(\l(-1)) \bar{t}(\l(0))^\psi \bar{t}(\l(1))^{\psi^2}\dots \bar{t}(\l(r))^{\psi^r}\dots$. Now the ring $A=A(G)$ is the subring of ${\rm rep}(G)$ spanned by all elements $t((\l))$, with $(\l) \in X^+(T)_\infty$. For $n \in \N_{-1}$, we write $A_n=A(G)_n$ for the subring of ${\rm rep}(G)$ spanned by all elements $t((\l))$ with $(\l) \in X^+(T)_n$ and $A=\cup_{n=1}^\infty A_n$. \\
	
	Let $X_s^+(T)$ be the set of all elements $\l \in X^+(T)$ such that either $(\l, \check{\al}) < \ell-1$ for some simple root $\alpha$ or $\ell-1 \leq (\l, \check{\beta}) \leq 2(\ell-1)$ for all simple roots $\beta$ of the quantum group $G$, where $\check{\al}=2\al/{(\al,\al)}$ and $\check{\beta}=2\beta/{(\beta,\beta)}$. Similarly let $X_s^+(\bar{T})$ be the set of all elements $\l \in X^+(T)$ such that either$(\l, \check{\al}) < p-1$ for some simple root $\al$ of $\bar{G}$ or $p-1 \leq (\l, \check{\beta}) \leq 2p-2$ for all simple roots $\beta$ of $\bar{G}$. We say that an element $(\l)=(\l(-1),\l(0),\dots,\l(r), \dots) \in X^+(T)_\infty$ is $n$-special if $\l(-1) \in X_s^+(T)$ and $\l(i) \in X_s^+(\bar{T})$ for all $0 \leq i \leq n$ and $\l(j)=0$ for all $j>n$. We let $X_s^+(T)_n$ for the subset of $X^+(T)_\infty$ consisting of all special $n$ elements. We say that an element $(\l) \in X^+(T)_\infty$ is special if it is $n$-special for sufficiently large $n$. We write $X_s^+(T)_\infty$ for the subset of $X^+(T)_\infty$ consisting of all special elements. \\
	
	For $\l \in X^+(T)$, we set ${\rm ln}(\l)= \sum(\l, \check{\beta})$ for all simple roots $\beta$ of the quantum group $G$. We call ${\rm ln}(\l)$ the length of $\l$. Note that ${\rm ln}(\l) \neq 0$ for all $0 \neq \l \in X^+(T)$. Similarly for $\l \in X^+(\bar{T})$, we have ${\rm ln}(\l)= \sum(\l, \check{\beta})$ where sum is over all the simple roots $\beta$ of $\bar{G}$. For each positive root $\al$ of quantum $GL_2$ (resp. for positive root $\al$ of $\bar{G}$), we have ${\rm ln}(\l)>0$ (resp. ${\rm ln}(\l)>0$) if $\l \neq 0$. We have the following proposition. \\ \\
	\textbf{Proposition 2.1} \\
\textit{(i) $A_n$ is spanned by $\{t((\l)) ~|~ (\l) \in X_s^+(T)_n\}$. \\
		(ii) $A$ is spanned by $\{t((\l)) ~|~ (\l) \in X_s^+(T)_\infty\}$.}

\begin{proof}
	We follow the proof of proposition in \cite{Dattag}. For an element $(\l)=(\l(-1),\l(0),\dots) \in X^+(T)_n$, we define ${\rm ln} ((\l))={\rm ln} (\l(-1))+ \ell (\Sigma_{i=0}^\infty p^i {\rm ln} (\l (i)))$. We prove the proposition by induction on the length of $(\l) \in X^+(T)_n$. For $n=-1, A_{-1}$ is spanned by all $t(\l), \l \in X_s^+(T)$ as tensor product of two tilting modules is again tilting. Assume now that $n>-1$. For $m \in \N_0$, let $A_{n,m}$ be the subgroup of $A_n$ generated by all elements $(\l) \in X^+(T)_n$ such that ${\rm ln} ((\l))=m$. If $m=0$, then $A_{n,m}=\Z[t(0)]$ and hence we are through. So assume the induction hypothesis upto $m-1$ and that $m>0$. \\
	
	Let $Z$ be the subgroup of $A_n$ spanned by all $t((\l)), (\l) \in X_s^+(T)_n$. We want to show that $Z=A_n$. Suppose that $(\l) \in X^+(T)_n$ such that ${\rm ln} ((\l))<m$. Then by induction we have $t((\l)) \in Z$. Hence we assume that $(\l)=(\l(-1),\l(0),\l(1),\dots) \in X^+(T)_n$ such that ${\rm ln} ((\l))=m$. \\ \\
	\textit{Case(i)}
	Suppose that $\l(-1) \in X_s^+(T)$. Let $(\mu)=(\l(0),\l(1),\dots)$. If $(\mu)=(0,0,\dots,0,\dots)$, then $T((\l))=T(\l(-1))$. Now by [\cite{Donqsc}, 72(40)], we have $t((\l))=t(\l(-1)) \in Z$. If $(\mu) \neq (0,0,\dots,0,\dots)$, then ${\rm ln} ((\mu))<m$. Hence by the induction we have $\bar{t}(\mu)=\Sigma_{(\xi)\in \Lambda}r_(\xi) \bar{t}(\xi)$, where $\Lambda$ is a finite subset of $X_s^+(T)_{n-1}$ consisting of elements of length at most ${\rm ln} ((\mu))$,and $r_{(\xi)} \in \Z$, for $(\xi) \in \Lambda$. Therefore we have $t((\l))=t(\l(-1))\bar{t}(\mu)^\phi$. Then $t((\l))=\Sigma_{(\xi)\in \Lambda'}r_(\xi) t((\xi^+))$, where $(\xi^+)=(\l(-1),\xi(0),\xi(1),\dots) \in X_s^+(T)_n$ and of length ${\rm ln}(\l(-1))+\ell {\rm ln} ((\xi)) \leq {\rm ln}(\l(-1))+\ell {\rm ln} ((\xi))={\rm ln} ((\l))=m$, as required. \\ \\
	\textit{Case(ii)}
	Suppose that $\l(-1) \notin X_s^+(T)$. Then we can write $\l(-1)=(\ell-1)\rho+\nu+\ell \tau$ with $\nu \in X_1(T)$ and $\tau \in X^+(\bar{T})$. Now consider the module ${\rm St} \otimes T(\nu) \otimes (T((\ell-1)\rho)+\ell \tau)$, which is a tilting module. By [\cite{Donqsc},p.no. 72], we have $T((\ell-1)\rho+\nu)$ is a direct summand of the tilting module ${\rm St} \otimes T(\nu)$ by and it occurs exactly once as a composition factor. Then again by [\cite{Donqsc},p.no. 73], we have $T((\ell-1)\rho+\nu) \otimes \bar{T}(\tau)^F$ is a tilting module with the highest weight $\l(-1)$. Hence we have $[T((\ell-1)\rho+\nu) \otimes T(\tau)^F]=[T(\l(-1))]+\Sigma_{(\xi)\in \Lambda}r_{(\xi)}[T(\xi)]$, where $\Lambda$ is a subset of $X^+(T)$ such that ${\rm ln} ((\xi))< {\rm ln} (\l(-1))$. This implies that $[T((\ell-1)\rho+\nu) \otimes \bar{T}(\tau)^F][\bar{T}(\mu)^F]=[T(\l(-1))]+\Sigma_{(\xi)\in \Lambda}r_{(\xi)}[T(\xi)][T(\mu)]^F$. Hence $[T(\l)]=[T((\ell-1)\rho+\nu) \otimes \bar{T}(\tau)^F][\bar{T}(\mu)^F]-\Sigma_{(\xi)\in \Lambda}r_{(\xi)}[T(\xi)][T(\mu)]^F$. For $(\xi) \in \Lambda$, let $(\eta)=(\xi,\l(0),\l(1),\dots)$. Now $(\eta) \in X^+(T)_n$ and
	\begin{equation*}
		{\rm ln}((\eta))={\rm ln} ((\xi))+{\rm ln}((\mu))<{\rm ln}(\l(-1))+\ell {\rm ln} ((\mu))= {\rm ln} ((\l))=m.
	\end{equation*}
		Hence by induction we have $t((\xi))t((\mu))^F=t((\eta))\in Z$. Therefore to show that $t((\l)) \in Z$, it is enough to show that $[T((\ell-1)\rho+\nu) \otimes \bar{T}(\tau)^F][\bar{T}(\mu)^F] \in Z$. Now by [\cite{Dattag}, section 1, Prop.], we have $t(\tau)t(\mu)=\Sigma_{(\theta) \in \Omega}s_\theta t((\theta))$, where $s_\theta \in \Z$ and $\Omega$ is a finite subset of $X_s^+(T)_n$ consisting of elements $(\theta)$ such that $\theta(-1)=0$. Hence we have 
		\begin{equation*}
		[T((\ell-1)\rho+\nu) \otimes \bar{T}(\tau)^F][\bar{T}(\mu)^F]= \Sigma_{(\theta )\in \Omega}s_\theta 	[T((\ell-1)\rho+\nu) \otimes \bar{T}(\theta)^F].
		\end{equation*}
	As each $[\bar{T}(\theta)] \in Z$ and $[T((\ell-1)\rho+\nu)] \in Z$ we have $[T((\ell-1)\rho+\nu) \otimes \bar{T}(\tau)^F][\bar{T}(\mu)^F] \in Z$ as required.
\end{proof}

\section{The structure of the tensor product of two simple modules of quantum $GL_2$}
In this section, we consider the quantum $GL_2$ at the primitive $\ell^{th}$ root of unity $q$ in $k$. Then $X(T)=\Z^2$ and $X^+(T)=\{ (\l_1,\l_2) ~|~ \l_1 \geq \l_2 \}$ and $W=S_2$. We let $X_1=\{ (a,b) \in X^+(T) ~|~ a-b \leq \ell -1 \}$ and $\bar{X_1}=\{ (a,b) \in X^+(\bar{T}) ~|~ a-b \leq 2p-2 \}$. Let $\pi= \{ (a,b) \in X^+(T) ~|~ a-b \leq 2(\ell -1) \}$ and $\bar{\pi}=\{ (a,b) \in X^+(\bar{T)} ~|~ a-b \leq 2p-2 \}$. Then $\pi$ and $\bar{\pi}$ are saturated subsets of $X^+(T)$ and $X^+(\bar{T})$ respectively. We write $\chi(\l)$ for the character of $\nabla(\l)$. Then we have the following results from \cite{Dattsm}: \\ \\
\textbf{Lemma 3.1} \ \textit{Let $\l=(\l_1,\l_2) \in \pi$. If $\l \in X_1$, then we have $L(\l)=T(\l)=\nabla(\l)=\Dt(\l)$. If $\l \in \pi \backslash X_1$, then we can write $\l_1-\l_2=\ell +r, \ 0 \leq r \leq \ell -2$. In this case $L(\l),L(\l-(r+1)(\eps_1-\eps_2))$ are the composition factors of $\nabla(\l)$. Further $\chi(\l)+\chi(\l-(r+1)(\eps_1-\eps_2))$ is the character of $T(\l)$.} \\ \\
\textbf{Lemma 3.2} \ \textit{Let $\l=(\l_1,\l_2), \l'=(\l_1',\l_2') \in X_1$. We can write $\l=(\l_1-\l_2)\eps_1+\l_2(\eps_1+\eps_2)$ and $\l'=(\l_1'-\l_2')\eps_1+\l_2'(\eps_1+\eps_2)$. We let $d_q=L(\eps_1+\eps_2), a=\l_1-\l_2, b=\l_1'-\l_2'$ and $r=\l_2+\l_2'$. Then $L(\l) \otimes L(\l')$ is $d_q^{\otimes r}$ times the direct sum of $T(u)$, where $u$ varies over the set $S(-1)$ which can be computed as follows: \\
	We list the weights $(a+b,0), (a+b-1,1),\dots,(a+b-i,i),\dots,(a+b-s,s)$ of $L(a,0)\otimes L(b,0)$, where $s=a+b-\ell, 0\leq s \leq (\ell -2)$. For each $j, 0 \leq j \leq s $, if $a+b-2j=\ell + u_j, 0 \leq u_j \leq (\ell -2)$ strike out $(a+b-j-u_j-1, j+u_j+1)$ from the list.} \\ \\
\textbf{Lemma 3.3} \ \textit{Let $\l=(\l_1,\l_2), \l'=(\l_1',\l_2') \in X^+(\bar{T})$ be such that $\l_1-\l_2, \l_1'-\l_2' < p$ and $d=L(\eps_1+\eps_2)$. Then $\bar{L}(\l) \otimes \bar{L}(\l')$ is $d^{\otimes(\l_2+\l_2')}$ times the direct sum of $\bar{T}(u)$, where $u$ varies over the set $S$ which can be computed as follows: Let $a=\l_1-\l_2, b=\l_1'-\l_2'$ and $a+b=p+r, 0 \leq r \leq (p-2)$. Let $A=\{ (a+b,0),(a+b-1,1),\dots,(a+b-r,r)\}$. Let $S(0)=A \backslash \{(a+b-i-u_i-1,i+u_i+1) ~|~ (a+b-i,i) \in A, a+b-2i = p+u_i, 0 \leq u_i \leq (p-2)\}$.} \\ \\ 
\textbf{Lemma 3.4} \ \textit{For $M \in {\rm mod}(G)$ and $N \in {\rm mod}(\bar{G})$, the map $\phi:M \otimes N^F \longrightarrow N^F \otimes M$ given by $\phi(x \otimes y)=y \otimes x$ is a $G$-module isomorphism.} \\ \\ 
The following theorem can be proved using lemma 3.2, lemma 3.3 and lemma 3.4. \\ \\
\textbf{Theorem 3.5} \ \cite{Dattsm} \textit{The tensor product of two simple $G$-modules is $d_q^{r'}(d^{r''})^F$ times a finite direct sum of indecomposable modules of the form $T(u_{-1}) \otimes (\bar{T}(u_0) \otimes \bar{T}(u_1)^{\bar{F}} \otimes \cdots \otimes \bar{T}(u_m)^{\bar{F}^m})^F$, where $u_{-1} \in X_1$ and $u_i \in \bar{\pi}$ for $0 \leq i \leq m$.} \\ \\
\textbf{Theorem 3.6} \ \cite{Dattsm} \textit{The indecomposable summand $T(u_{-1}) \otimes (\bar{T}(u_0) \otimes \bar{T}(u_1)^{\bar{F}} \otimes \cdots \otimes \bar{T}(u_m)^{\bar{F}^m})^F$ in the preceding theorem is tilting if and only if $\ell -1 \leq a_{-1}-b_{-1} \leq 2(\ell -1)$ and $p-1 \leq a_i-b_i \leq 2p-2$, where $u_i=(a_i,b_i)$ for $i=-1,0,\dots,m-1$. Also $T(u)$ is simple if $u_i \in \{ (a,b) \in X^+(\bar{T}) ~|~ a-b \leq p-1 \}, i=0,1,\dots,m.$} 

\section{Generators and Relations.}

Now we consider the ring $A$ of twisted tilting modules of the quantum $GL_2(k)$. We give a set of generators and relations for this ring $A$. As in Section 3, we identify 
 $X^+(T)_n,$ $ X^+(T)_{\infty},$ $ X^+_s(T)_n,$  $X^+_s(T)_{\infty}$ as subsets of $X^+(T)\times X^+(\b{T})_{\infty} $ in the usual way.  By [\cite{Donqsc}, p.no 72 (5)] and [\cite{Donqsc},3.4 (1), p. no. 75],  $\{T(\l_1,\l_2))~|~(\l_1,\l_2)\in X^+(T), \l_1-\l_2\leq 2l-2\}$ is a set of pairwise non-isomorphic indecomposable $G_1$-modules. \\

\ni {\bf Proposition 4.1}

\ni {\it (i) For $a\in X^+_s(T)_n$, the twisted tilting module $T(a)$ is indecomposable.\\
(ii)  For $a,b\in X^+_s(T)_n$ the twisted tilting modules $T(a)$ and $T(b)$ are isomorphic if and only $a=b$.\\
(iii) Any indecomposable twisted tilting modules of height at most $n$ is isomorphic to $T(a)$, for some $a\in X^+_s(T)_n$.\\ 
(iv) The ring $A_n$ is a free ${\mathbb{Z}}$-module with ${\mathbb{Z}} $ basis $\{T(a)~|~a\in X^+_s(T)_n\}$.\\
(v) The ${\mathbb{Z}}$-module $A_n$ is free of rank $(2l-2)(2p-2)^n$ as a module over $\psi ^{n+1}(\b{A}_0)$.\\
(vi) The ring $A$ is a free ${\mathbb{Z}}$-module with ${\mathbb{Z}}$-basis $\{t(a)~|~a\in X^+_s(T)_{\infty}\}$.}\\

\begin{proof}
(i) Let $a\in X^+_s(T)_n$. If $n=-1$, then $T(a)$ is an indecomposable tilting module. So, let us assume that $n\geq 0$. Let $\b{a}=(a(0),a(1),\dots ,a(n))\in X^+_s(\b{T})_n$ and $a(i)=(a(i)_1,a(i)_2)\in {\mathbb{Z}}^2$, $i\geq -1$. Then by the [\cite{Dattag}, Section 2, Prop], we have $\b{T}(\b{a})$ is indecomposable. As $0\leq a(-1)_1-a(-1)_2\leq 2(l-1)$, by [\cite{Donqsc}, p.no.72,(5)], we have $T(a)=T(a(-1))\ot (\b{T}(\b{a}))^F$ is indecomposable.\\

\ni (ii) Let $a,b\in X^+_s(T)_n$. Suppose that $T(a)\cong T(b)$. We want to show that $a=b$. Let $a=(a(-1),a(0),\cdots ,a(n))$ and $b=(b(-1),b(0),\cdots ,b(n))$. As $G_1$-module, we have $T(a(-1))=T(a)_{|G_1}\cong T(b)_{|G_1}= T(b(-1))$. The highest weight $a(-1)$ of $T(a(-1))$ is same as the highest weight $b(-1)$ of $T(b(-1))$. Let $L={\rm soc}(T(a(-1))={\rm soc}(T(b(-1))$, which is simple. Now we have $${\rm Hom}_{G_1}(L,T(a(-1))\ot \b{T}(\b{a})^F)\cong {\rm Hom}_{G_1}(L,T(b(-1))\ot \b{T}(\b{b})^F )$$
$${\rm End}_{G_1}(L)\ot \b{T}(\b{a})^F\cong {\rm End}_{G_1}(L)\ot \b{T}(\b{b})^F$$

Hence we have $$\b{T}(\b{a}))\cong \b{T}(\b{b})$$

Again by [\cite{Dattag}, Section 2, Prop], we have $\b{a}=\b{b}$. Therefore $a=b$.\\

\ni (iii), (iv), (v) and (vi) follows from Proposition of Section 1 of \cite{Dattag}.
\end{proof} 

Let $J={\mathbb{Z}}[d_q,d]$. For $r\in {\mathbb{N}}_0$, we have Laurent polynomials $P_r(t)\in J[t,t^{-1}]$ defined recursively as follows: $P_0(t)=1$, $P_1(t)=t+t^{-1}$ and $$(t+t^{-1})P_{r}(t)=P
_{r+1}(t)+d_qP_{r-1}(t)$$

Similarly, we also consider the following polynomials $Q_r(s), r\in {\mathbb{N}}_0$, recursively in the ring $J[s,s^{-1}]$ as follows: 
We let $Q_0=1$, $Q_1=s+s^{-1}$ and for $r\geq 2$, $$(s+s^{-1})Q_r=Q_{r+1}(s)+dQ_{r-1}(s)$$ 

Let $\omega =t+t^{-1}$ and $\upsilon=s+s^{-1}$. If a Laurent polynomial $\sum _{i=-\infty}^{\infty}a_it^i\in J[t,t^{-1}]$ (resp. $J[s,s^{-1}]$) is such that the coefficient of $t^i$ is same as the coefficient of $t^{-i}$  for all $i$, then it can be expressed as a polynomial in $J[\omega]$ (resp. as a polynomial in $J[\upsilon]$). \\

Let $P_r(\omega)=P_r(t)$ and $Q_r(\upsilon)=Q_r(s)$ for $r\geq 0$. Then we have $P_0(\omega)=1$, $P_1(\omega)=\omega$ and $\omega P_r(\omega)=P_{r+1}(\omega)+d_qP_{r-1}(\omega)$. 
Similarly we have $Q_0(\upsilon)=1$, $Q_1(\upsilon)=\upsilon$ and for $r\geq 2$, we have $\upsilon Q_r(\upsilon)=Q_{r+1}(\upsilon)+dQ_{r-1}(\upsilon)$. We let $g(t)=t^l+t^{-l}$ and $h(s)=s^p+s^{-p}$. Now we have $g'(\omega)=lP_{l-1}(\omega)$ and $h'(\upsilon)=pQ_{p-1}(\upsilon)$ and also $P_{2l-1}(\omega)=g(\omega)P_{l-1}(\omega)$ and $Q_{2p-1}(\upsilon)=h(\upsilon)Q_{p-1}(\upsilon)$.\\

For any Laurent polynomial $f(t)\in J[t,t^{-1}]$, we simply write $f(E)$ for the element $f([E])$ of $A$. Now we identify the integral group ring ${\mathbb{Z}}X(T)$ with  a subring of $J[t,t^{-1}]$ by $\epsilon _1\mapsto t$ and $\eps_1+\eps_2\mapsto d_q$. Now for any $(r,0)\in X^+(T)$, we have ${\rm ch}\nabla(r,0)=P_r(t)$, $r\geq 0$.  Therefore the element $P_r(E)$ has same character as the $G$-module $\nabla(r,0)$. By [\cite{Donqsc},(7), p.no.73], we know that $$T(2l-1,0)\cong T((l-1,0)\otimes \b{E}^F$$ And similarly by the [\cite{Donqsc}, 3.4(1), p.no. 75], we have $\b{T}(2p-1,0)\cong \b{T}(p-1,0)\otimes \b{E}^{\b{F}}$ and hence we have $\b{T}(2p-1,0)^F\cong \b{T}(p-1,0)^F\otimes (\b{E}^{\b{F}})^F$.\\

Let $\wt{J}={\mathbb{Z}}[d_q,d,d_{q}^{-1},d^{-1}]$ and $K$ be its field of fractions. Here we consider the free polynomial ring $\wt{J}[X_{-1},X_0,X_1,\cdots ]$. We define a map $\Phi:\wt{J}[X_{-1},X_0,\dots ]\lr A$ by $\Phi(X_{-1})=[E]$, for $i\geq 0$, $\Phi(X_i)=[(\b{E}^{\b{F}})^F]$ and $\Phi(d_q)=D_q, \Phi(d)=D$. Now by the above discussion, we have $D_q^l-D^l, P_{l-1}(X_0-g(X_{-1}))$ and $Q_{p-1}(X_i)(X_{i+1}-h(X_i))\in \Ker \Phi$ for all $i\geq 0$. Let $\Phi _n$ be the restriction of $\Phi$ to the ring $\wt{J}[X_{-1},X_0,\dots ,X_n]$, then by [Prop.4.1 ( iv)], we have the following. \\

\ni {\bf Proposition 4.2.} {\it For $n\geq -1$, the restriction $\Phi_{n}:\wt{J}[X_{-1},X_0,\dots, X_n]\lr A_n$ be the restriction of the ring homomorphism. Then $\Phi_n$ is surjective and  $\Ker (\Phi_n)$ is  generated by the elements  $(P_{l-1}(X_{-1})(X_0-g(X_{-1}))$ and $(Q_{p-1}(X_i)(X_{i+1}-h(X_i)))$, for $0\leq i<n$.} \hfill{$\Box$}\\ 

Now we show that the ring is reduced. Let $\Phi(X_i)=x_i$, $i\geq -1$. Let ${\mathbb{C}}$ be the field of complex numbers. We let  $A_{\mathbb{C}}(n)={\mathbb{C}}\otimes _{\mathbb{Z}}A_n$.  By \cite{Dattag}, we know that the polynomials $Q_{p-1}(X_{i+1}-h(X_i)))$, $i\geq 0$, has no repeated roots. Hence to show that the ring $A_{\mathbb{C}}(n)$ is reduced, it is enough to show that $(P_{l-1}(X_{-1})(X_0-g(X_{-1}))$ has no repeated roots. \\

First we note the ring $A_{-1}=\wt{J}[x_{-1}]$ is a reduced ring. Any element $F\in A_{\mathbb{C}}(0)$ is uniquely expressible in the form $F=\sum _{i=0}^{2l-2}a_iP_i$, $a_i\in J[x_0,x_1,\cdots ,x_n]$.\\

\ni {\bf Theorem 4.3} {\it The ring $A$ is reduced.}\\

\ni The theorem follows from the following Lemma.\\

\ni {\bf Lemma 4.4} {\it For any $b\in {\mathbb{Z}}$, $(x_{-1}-b)$ is not a zero-divisor in $A_K(0)$}.

\begin{proof}
Suppose $b\in {\mathbb{Z}}$ and $(x_{-1}-b)F=0$. Then we have $$x_{-1}(\sum_ {i=0}^{2l-2}a_iP_i)=bF$$
\ni Hence $$\sum _{i=0}^{2l-2}a_ix_{-1}P_i=bF$$ We know that $x_{-1}P_i=D_qP_{i-1}+P_{i+1}$, for $i\geq 2$

\ni Hence we have $$a_0x_{-1}P_0+\sum _{i=1}^{2l-3}a_ix_{-1}P_i+a_{2l-2}x_{-1}P_{2l-2}=bF$$

\ni It follows that $$a_0P_1+\sum _{i=1}^{2l-3}a_i(P_{i+1}+dP_{i-1})+a_{2l-2}(P_{2l-1}+dP_{2l-3})=bF$$

\ni Therefore we have $$a_0P_1+\sum_{i=2}^{2l-2}a_{i-1}P_i+d\sum _{i=0}^{2l-4}a_{i+1}dP_i+a_{2l-2}(P_{2l-1}+dP_{2l-3})=bF$$

\ni And hence we have $$a_0P_1+\sum _{i=2}^{2l-2}a_{i-1}P_i+d\sum _{i=0}^{2l-4}a_{i+1}dP_i+a_{2l-2}x_0P_{l-1}+da_{2l-2}P_{2l-3}=bF$$

\ni By comparing the coefficients, we have\\

$$d^ia_i=P_i(b)a_0 \hspace{0.75cm}\mbox{for}\hspace{0.5cm}0\leq i\leq l-1\hspace{2cm}(\star)\\$$
$$a_{2l-2-i}=P_i(b)a_{2l-2}\hspace{0.4cm}\mbox{for}\hspace{0.5cm}0\leq i\leq l-1\hspace{1.5cm}(\star\star)$$

$$a_{l-2}+da_l+a_{2l-2}x_0=ba_{l-1}$$

\ni Now by $(\star)$ and $(\star\star)$, we have \\

$$P_{l-2}(b)a_{2l-2}+x_0a_{2l-2}=bP_i(b)a_{2l-2}.$$

\ni Hence we have $$(x_0-(P_l(b)-P_{l-2}(b))\displaystyle{\frac{a_0}{d^{l-1}}}=0.$$

\ni Therefore $$ (x_0-(P_l(b)-P_{l-2}(b))a_0=0.$$

If $a_0=0$, then $a_{2l-2}=0$ and hence $F=0$. If $a_0\neq 0$, then $(x_0-(P_l(b)-P_{l-2}(b))$ is a zero divisor in $A_K(0)$. Now by [\cite{Dattag}, Section 3, Prop.], for $d=1$ and $q=1$ and $b=2$, we have $(x_0-2)$ is a zero divisor which is not true. Therefore $x_0\pm 2d$ are not zero divisor in $A_K(0)$.
\end{proof}

\begin{center}
 {\bf Statement of Conflict of Intereest}
\end{center}

 The authors do not have any conflict of interest with any of the content of this article.\\

\end{document}